\documentclass{ws-procs975x65}            
\begin{document}
\title{Chern-Simons Theory, Vassiliev Invariants, Loop Quantum Gravity and Functional Integration Without Integration}

\author{Louis H. Kauffman}

\address{Mathematics Department, University of Illinois at Chicago,\\
Chicago, IL 60607-7045, USA\\
$^*$E-mail: kauffman@uic.edu\\
www.math.uic.edu/~kauffman}

\begin{abstract}This paper is an exposition of the
relationship between  Witten's Chern-Simons functional integral and the theory of Vassiliev invariants of knots and links
in three dimensional space. We conceptualize the functional integral in terms of equivalence classes
of functionals of gauge fields and we do not use measure theory. This approach makes it possible to discuss the mathematics intrinsic to the functional integral rigorously and without 
functional integration. Applications to loop quantum gravity are discussed.
\end{abstract}


\bodymatter
\section{Introduction}

\noindent This paper is an introduction to how  Vassiliev
invariants in knot theory arise naturally in the context of
Witten's functional integral. The relationship between Vassiliev invariants and
Witten's integral has been known since Bar-Natan's thesis \cite{Bar-Natan-Thesis}
where he discovered, through this connection, how to define Lie algebraic weight
systems for these invariants. \\

This paper is written in a context of ``integration without integration''. 
The idea is as follows. Let $F(A),G(A),H(A)$ be functionals of a gauge field $A$ that vanish rapidly as the amplitude 
of the field goes to infinity. We say that $F ~\sim~ G$ if $F - G = DH$ where $D$ denotes a gauge functional derivative.
We define $\int F(A)$ to be the equivalence class of $F(A).$ By definition, this integral satisfies integration by parts,
and it is a useful conceptual substitute for a functional integral over all gauge fields (modulo gauge equivalence).
We replace the usual notion of functional integral with such equivalence classes. \\

The paper is a sequel to  \cite{Heuristics} and \cite{WittKont}. 
In these papers we show somewhat more about the relationship of Vassiliev invariants
and the Witten functional integral. In particular, we show how the Kontsevich integrals (used to
to give rigorous definitions of these invariants) arise as Feynman integrals in the 
perturbative expansion of the Witten functional integral.
See also the work of Labastida and P$\acute{e}$rez \cite{LP} on this same subject. 
The result is an interpretation of the Kontsevich integrals in terms of the light-cone
gauge and thereby extending the original work of Fr\"ohlich and King
\cite{Frohlich and King}.  The purpose of this paper is to give an exposition of
the beginnings of these relationships, to introduce  diagrammatic techniques that illuminate the
connections, and to show how the integral can be fruitfully formulated in terms of certain equivalence
classes of functionals of gauge fields. \\

The paper is divided into six sections beyond the introduction.   
Section 2 discusses Vassiliev invariants and invariants of rigid vertex graphs.  
Section 3 discusses the concept of replacing integrals by equivalence classes. Section 4 
introduces the basic formalism and shows how the functional integral, regarded without integration, is related directly to knot invariants and particularly,
Vassiliev invariants. Section 5 discusses the formalism of the perturbative expansion of the Witten integral. 
Section 6 is a sketch of the loop transform, useful in loop quantum gravity and ends with a quick discussion of the Kodama state with references to recent literature. 
Section 7 discusses how the Kontsevich integrals for Vassiliev invariants arise from the 
perturbation expansion. 
 \\

\noindent {\bf Acknowledgement.}  We thank students and colleagues for many stimulating conversations on the themes of this paper, and we thank the organizers of the
Conference on 60 Years of Yang-Mills Gauge Field Theories (25 to 28 May 2015) for the invitation and opportunity to speak about these ideas in Singapore.\\

\bigbreak

\section{Vassiliev Invariants and Invariants of Rigid Vertex Graphs}

If  $V(K)$ is a  (Laurent polynomial valued,  or more generally - commutative
ring valued) invariant of knots,  then it can be naturally extended to an
invariant of rigid vertex graphs \cite{Kauffman-Graph} by defining the invariant
of graphs in terms of the knot invariant via an `unfoldingÓ  of the vertex. That
is, we can regard the vertex as a `black box" and replace it by any tangle of our
choice. Rigid vertex motions of the graph preserve the contents of the black box,
and hence implicate ambient isotopies of the link obtained by replacing the black
box by its contents. Invariants of knots and links that are evaluated on these
replacements are then automatically rigid vertex invariants of the corresponding
graphs. If we set up a collection of multiple replacements at the vertices with
standard conventions for the insertions of the tangles, then a summation over all
possible replacements can lead to a graph invariant with new coefficients
corresponding to the different replacements.  In this way each invariant of knots
and links implicates a large collection of graph invariants. See
\cite{Kauffman-Graph}, \cite{Kauffman-Vogel}. \vspace{3mm}

The simplest tangle replacements for a 4-valent vertex are the two crossings,
positive and negative, and the oriented smoothing. Let V(K) be any invariant of
knots and links. Extend V to the category of rigid vertex embeddings of 4-valent
graphs by the formula $$V(K_{*}) = aV(K_{+}) + bV(K_{-}) + cV(K_{0})$$ where
$K_{+}$ denotes a knot diagram $K$ with a specific choice of positive crossing,
$K_{-}$ denotes a diagram identical to the first with the positive crossing
replaced by a negative crossing and  $K_{*}$ denotes a diagram identical to the
first with the positive crossing replaced by a graphical node. \vspace{3mm}

This formula means that we define  $V(G)$  for an embedded 4-valent graph  $G$ 
by taking the sum

$$V(G) = \sum_{S} a^{i_{+}(S)}b^{i_{-}(S)}c^{i_{0}(S)}V(S)$$

\noindent with the summation over  all knots and links $S$ obtained from  $G$ by
replacing a node of $G$ with either a crossing of positive or negative type, or
with  a smoothing of the crossing that replaces it by a planar embedding of
non-touching segments (denoted $0$).  It is not hard to see that if $V(K)$  is an
 ambient isotopy invariant of knots, then,  this extension is an rigid vertex
isotopy invariant of graphs.  In rigid vertex isotopy the cyclic order at the
vertex is preserved, so that the vertex behaves like a rigid disk with flexible
strings attached to it at specific points. \vspace{3mm}

There is a rich class of graph invariants that can be studied in this manner. 
The Vassiliev Invariants \cite{Birman and Lin},\cite{Bar-Natan}
constitute the important special case of these graph invariants where  $a=+1$,
$b=-1$ and $c=0.$    Thus  $V(G)$  is a Vassiliev invariant if

$$V(K_{*}) = V(K_{+})  -  V(K_{-}).$$

\noindent Call this formula the {\em exchange identity} for the Vassiliev
invariant $V.$  See Figure 1 \vspace{3mm}

\begin{figure}[htb]
     \begin{center}
     \begin{tabular}{c}
     \includegraphics[width=6cm]{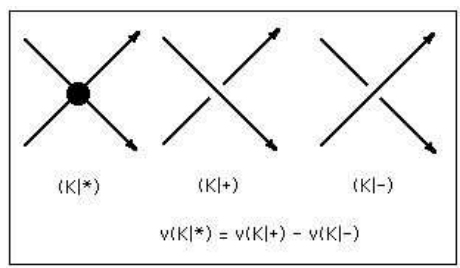}
     \end{tabular}
     \caption{\bf Exchange Identity for Vassiliev Invariants}
     \label{Figure 1}
\end{center}
\end{figure}

$V$  is said to be of  {\em finite type}  $k$  if  $V(G) = 0$  whenever  $|G| >k$
 where $|G|$  denotes the number of (4-valent) nodes in the graph $G.$ The notion
of finite type is of extraordinary significance in studying these invariants. One
reason for this is the following basic Lemma. \vspace{3mm}

\noindent {\bf Lemma.} If a graph $G$ has exactly $k$ nodes, then the value of a
Vassiliev invariant $v_{k}$ of type $k$ on $G$, $v_{k}(G)$, is independent of the
embedding of $G$. \vspace{3mm}

\noindent {\bf Proof.} The different embeddings of $G$ can be represented by link
diagrams with some of the 4-valent vertices in the diagram corresponding to the
nodes of $G$. It suffices to show that the value of $v_{k}(G)$ is unchanged under
switching of a crossing.  However, the exchange identity for $v_{k}$ shows that
this difference is equal to the evaluation of $v_{k}$ on a graph with $k+1$ nodes
and hence is equal to zero. This completes the proof.// \vspace{3mm}

The upshot of this Lemma is that Vassiliev invariants of type $k$ are intimately
involved with certain abstract evaluations of graphs with $k$ nodes. In fact,
there are  restrictions (the four-term relations) on these evaluations demanded
by the topology  and it follows from results of Kontsevich \cite{Bar-Natan} that
such abstract evaluations actually determine the invariants. The knot  invariants
derived from classical Lie algebras are all built from Vassiliev invariants of
finite type. All of this is directly related to Witten's functional integral
\cite{Witten}. \vspace{3mm}

In the next few figures we illustrate some of these main points. In Figure 2 we
show how one associates a so-called chord diagram to represent the abstract graph
associated with an embedded graph. The chord diagram is a circle with arcs
connecting those points on the circle that are welded to form the corresponding
graph.  In Figure 3 we illustrate how the four-term relation is a consequence of
topological invariance. In Figure 4 we show how the four term relation is a
consequence of the abstract pattern of the commutator identity for a matrix Lie
algebra. This shows that the four term relation is directly related to a
categorical generalisation of Lie algebras. Figure 5 illustrates how the weights
are assigned to the chord diagrams in the Lie algebra case -  by inserting Lie
algebra matrices into the circle and taking a trace of a sum of matrix products.
\vspace{3mm}

\begin{figure}[htb]
     \begin{center}
     \begin{tabular}{c}
     \includegraphics[width=6cm]{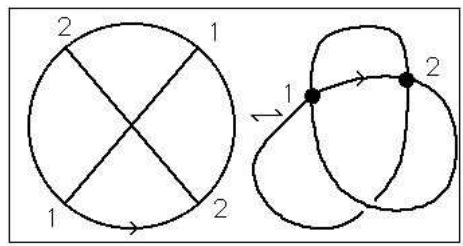}
     \end{tabular}
     \caption{\bf Chord Diagrams}
     \label{Figure 2}
\end{center}
\end{figure}

\begin{figure}[htb]
     \begin{center}
     \begin{tabular}{c}
     \includegraphics[width=6cm]{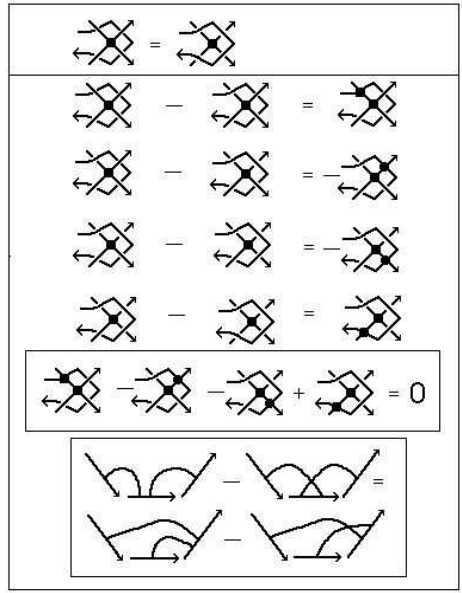}
     \end{tabular}
     \caption{\bf The Four Term Relation from Topology}
     \label{Figure 3}
\end{center}
\end{figure}

\begin{figure}[htb]
     \begin{center}
     \begin{tabular}{c}
     \includegraphics[width=6cm]{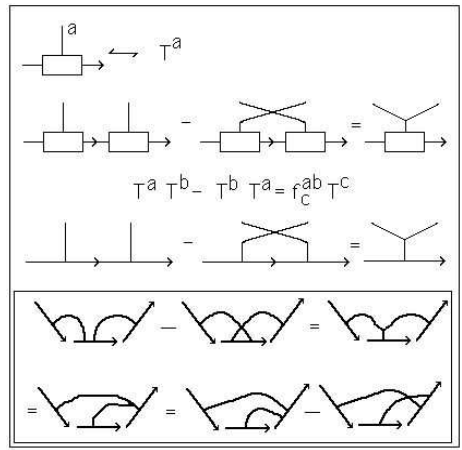}
     \end{tabular}
     \caption{\bf The Four Term Relation from Categorical Lie
Algebra}
     \label{Figure 4}
\end{center}
\end{figure}

\begin{figure}[htb]
     \begin{center}
     \begin{tabular}{c}
     \includegraphics[width=6cm]{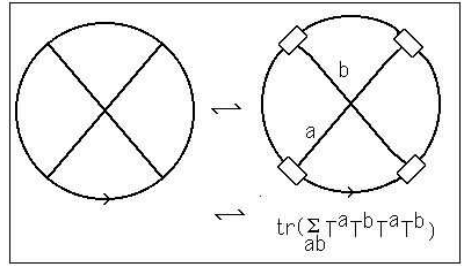}
     \end{tabular}
     \caption{\bf Calculating Lie Algebra Weights}
     \label{Figure 5}
\end{center}
\end{figure}

\section{Integration without integration}
Recall that if $Z = \int_{-\infty}^{\infty} e^{-x^{2}/2}dx$ then
$$Z^2 = \int_{-\infty}^{\infty}\int_{-\infty}^{\infty} e^{-(x^{2}+y^{2})/2}dxdy
=\int_{0}^{2\pi}\int_{0}^{\infty} e^{-r^{2}/2}rdrd\theta$$
$$=2\pi \int_{0}^{\infty} e^{-r^{2}/2}rdr = 2\pi.$$ Whence $$Z = \sqrt{2 \pi}.$$ Furthermore,
if $$Z(J) = \int_{-\infty}^{\infty} e^{-x^{2}/2 + Jx}dx,$$ then 
$$Z(J) = \int_{-\infty}^{\infty} e^{-(x-J)^{2}/2 + J^{2}/2}dx$$
$$=e^{J^{2}/2}\int_{-\infty}^{\infty} e^{-(x-J)^{2}/2}dx$$
$$=e^{J^{2}/2}\int_{-\infty}^{\infty} e^{-x^{2}/2}dx$$
$$=e^{J^{2}/2}Z(0) = \sqrt{2\pi} e^{J^{2}/2}.$$  Now examine how much of this calculation could be done if we did
not know about the existence of the integral, or if we did not know how to calculate explicitly the values of these 
integrals across the entire real line. Given that we believed in the existence of the integrals, and that we
could use properties such as change of variable giving
$$\int_{-\infty}^{\infty} e^{-(x-J)^{2}/2}dx = \int_{-\infty}^{\infty} e^{-x^{2}/2}dx,$$ we could deduce the
relative result stating that 
$$Z(J) = \int_{-\infty}^{\infty} e^{-x^{2}/2 + Jx}dx = e^{J^{2}/2}\int_{-\infty}^{\infty} e^{-x^{2}/2}dx.$$ From this we
can deduce that 
$$d^{n}Z(J)/dJ^{n}|_{J=0} = d^{n}/dJ^{n}\int_{-\infty}^{\infty} e^{-x^{2}/2 + Jx}dx = \int_{-\infty}^{\infty}
x^{n}e^{-x^{2}/2}dx.$$  Hence
$$\int_{-\infty}^{\infty}x^{n}e^{-x^{2}/2}dx = d^{n}(e^{J^{2}/2})/dJ^{n}|_{J=0}\int_{-\infty}^{\infty} e^{-x^{2}/2}dx.$$

But now, lets go a step further and imagine that we really have no theory of integration available. Then we are 
in the position of freshman calculus where one defines $\int f$ to be ``any" function $g$ such that $dg/dx = f.$
One {\it defines} the integral in this form of elementary calculus to be the anti-derivative, and this takes care
of the matter for a while! What are we really doing in freshman calculus? We are noting that for integration
on an interval $[a,b],$ if two functions $f$ and $g$ satisfy 
$f - g = dh/dx$ for some differentiable function $h,$ then we have that 
$$\int_{a}^{b}(f-g) = \int_{a}^{b}dh/dx = h(b)-h(a).$$  If the function $h(x)$ vanishes as $x$ goes to infinity, then
we have that $$\int_{-\infty}^{\infty} fdx = \int_{-\infty}^{\infty} gdx$$ when $f-g = dh/dx.$ This suggests turning
things upside down and {\it defining an equivalence relation on functions} $$f \sim g$$ if $$f-g = dh/dx$$ where $h(x)$
is a function vanishing at infinity. Then {\it we define the integral  $$\int f(x)$$ to be the equivalence class of the
function $f(x).$} This ``integral" represents integration from minus infinity to plus infinity but it is defined only as
an equivalence class of functions. An ``actual" integral, like the Riemann, Lesbeque or Henstock integral is a
well-defined real valued function that is constant on these equivalence classes.
\bigbreak

We shall say that 
$f(x)$ is {\it rapidly vanishing at infinity} if $f(x)$ and all its derivatives are vanishing at infinity.
For simplicity, we shall assume that all functions under consideration have convergent power series expansions so
that $f(x + J) = f(x) + f'(x)J + f''(x)J^{2}/2! + \cdots,$ and that they are rapidly vanishing at infinity. It then
follows that 
$$f(x+J) = f(x) + d(f(x)J + f'(x)J^{2}/2! + \cdots)/dx \sim f(x),$$ and hence we have that 
$\int f(x + J) = \int f(x),$ giving translation invariance when $J$ is a constant. 
\bigbreak

\noindent We have shown the following Proposition.
\smallbreak

\noindent {\bf Proposition.} Let $f(x), g(x), h(x)$ be functions rapidly vanishing at infinity (with power series
representations). Let $\int f$ denote the equivalence class of the function $f$ where $ f \sim g$ means that 
$f-g = Dh$ where $Dh = dh/dx.$ Then this integral satisfies the following properties
\begin{enumerate}

\item If $f \sim g$ then $\int f = \int g$.

\item If $k$ is a constant, then $\int (kf+g) = k\int f + \int g.$ 

\item  If $J$ is a constant, then $\int f(x + J) = \int f(x).$

\item $\int Dh = 0$ where $0$ denotes the equivalence class of the zero function. Hence
$\int f(Dg) + \int (Df)g = \int D(fg) = 0,$ so that integration by parts is valid with vanishing boundary conditions
at infinity.

\end{enumerate}
\smallbreak

Note that $ e^{-x^{2}/2}$ is rapidly vanishing at infinity.  
We now see that most of the calculations that we made about $e^{-x^{2}/2}$ were actually statements about the 
equivalence class of this function:
$$e^{-x^{2}/2 + Jx} = e^{-(x-J)^{2}/2 + J^{2}/2} = e^{J^{2}/2}e^{-(x-J)^{2}/2} \sim e^{J^{2}/2}e^{-x^{2}/2},$$ whence
$$\int e^{-x^{2}/2 + Jx} =  e^{J^{2}/2} \int e^{-x^{2}/2}.$$

\subsection{Functional Derivatives}
In order to generalize the ideas presented in this section to the context of functional integrals, we need to 
discuss the concept of functional derivatives. We are given a {\it functional} $F(\alpha(x))$ whose argument $\alpha(x)$
is a function of a variable $x.$ We wish to define the {\it functional derivative} $\delta F(\alpha(x))/\delta
\alpha(x_0)$ of $F(\alpha(x))$ with respect to $\alpha(x)$ at a given point $x_{0}.$ The idea is to regard each
$\alpha(x_{0})$ as a separate variable, giving
$F(\alpha(x))$ the appearance of a function of infinitely many variables. In order to formalize this notion one needs to
use  generalized functions (distributions) such as the Dirac delta function $\delta(x)$, a distribution with the
property that $\int_{a}^{b}\delta(x_{0})f(x)dx = f(x_{0})$ for any integrable function $f(x)$ and point $x_0$ in the
interval $[a,b].$ One defines the functional derivative by the formula

$$\delta F(\alpha(x))/\delta \alpha(x_0) = lim_{\epsilon \rightarrow 0}[F(\alpha(x) + \delta(x_0)\epsilon) -
F(\alpha(x))]/\epsilon.$$

\noindent Note that if 

$$F(\alpha(x)) = \alpha(x)^2$$ then 

$$\delta F(\alpha(x))/\delta \alpha(x_0) = lim_{\epsilon \rightarrow 0}[(\alpha(x) + \delta(x_0)\epsilon)^2 -
\alpha(x)^2]/\epsilon$$
$$= lim_{\epsilon \rightarrow 0}[2\alpha(x)\delta(x_0)\epsilon + \delta(x_0)^{2}\epsilon^{2}]/\epsilon$$
$$= 2\alpha(x)\delta(x_0).$$ While if

$$G(\alpha(x)) = \int_{a}^{b}\alpha(x)^{2}dx$$ then 
$$\delta G(\alpha(x))/\delta \alpha(x_0) = 2\alpha(x_0)$$ when $x_0 \in [a,b].$ More generally, if
$$G(\alpha(x)) = \int_{a}^{b}f(\alpha(x))dx$$ for a differentiable function $f,$ then
$$\delta G(\alpha(x))/\delta \alpha(x_0) = f'(\alpha(x_{0})).$$ These examples show that the results of a functional
differentiation can be either a distribution or a function, depending upon the context of the original
functional.
\bigbreak

In the case of a path integral of the type used in quantum mechanics, one wants to integrate a functional $F(p)$
over paths $p(t)$ with $t$ in an interval $[0,1].$ The functional takes the form
$$F(p) = e^{(i/\hbar)\int_{0}^{1}S(p(t))dt}$$ and the traditional Feynman path integral has the form
$$\int dP e^{(i/\hbar)\int_{0}^{1}S(p(t))dt},$$ giving the amplitude for a particle to travel from $a = p(0)$ to
$b = p(1),$ the integration proceeding over all paths with these initial and ending points. 

Here the equivalence relation corresponding to the functional integral is $F \sim G$ if $F - G = DH$ where
$$DH = \delta H(p)/\delta p(t_0)$$ for some time $t_{0}$ and some $H(p).$ Again we need to specify the 
class of functionals and to say what it means for a functional to "vanish at infinity." Since we are integrating over
all paths, we need a notion of size for a path. This can be defined by $$||p|| = (\int_{0}^{1}|p(t)|^{2}dt)^{1/2}.$$
Note that for $F(p) = e^{(i/\hbar)\int_{0}^{1}S(p(t))dt}$ we have
$$\delta F(p)/\delta p(t_0) = (i/\hbar)[\delta\int_{0}^{1}S(p(t))dt/\delta p(t_0)] F(p)$$
$$= (i/\hbar)S'(p(t_0))F(p).$$ Here we see the fact that the integral can be dominated by contributions from paths
where this variation is zero. Note that in order to estimate this stationary phase contribution to the functional
integral, one needs more than just a definition of the integral as an equivalence class of functionals. Nevertheless, 
we shall see in the next section that these equivalence classes do give insight into the topology associated with
Witten's integral.

\section{Vassiliev Invariants and Witten's Functional Integral}

In \cite{Witten}  Edward Witten proposed a formulation of a class of 3-manifold
invariants as generalized Feynman integrals taking the form  $Z(M)$  where

$$Z(M) = \int DAe^{(ik/4\pi)S(M,A)}.$$

\noindent Here  $M$ denotes a 3-manifold without boundary and $A$ is a gauge
field  (also called a gauge potential or gauge connection)  defined on $M$.  The
gauge field is a one-form on a trivial $G$-bundle over  $M$ with values in a
representation of the  Lie algebra of $G.$ The group $G$ corresponding to this
Lie algebra is said to be the gauge group. In this integral the ÒactionÓ  
$S(M,A)$  is taken to be the integral over $M$ of the trace of the Chern-Simons
three-form   $A \wedge dA + (2/3)A \wedge A \wedge A$.  (The product is the wedge
product of differential forms.) \vspace{3mm}

$Z(M)$  integrates over all gauge fields modulo gauge equivalence.
\vspace{3mm}

The formalism  and   internal logic of Witten's integral supports  the existence
of a large class of topological invariants of 3-manifolds and  associated
invariants of knots and links in these manifolds. \vspace{3mm}

The invariants associated with this integral have been given rigorous
combinatorial descriptions but questions and conjectures arising from the integral
formulation are still outstanding. Specific conjectures about this integral take the form
of just how it implicates invariants of links and 3-manifolds, and how these invariants
behave in certain limits of the coupling constant $k$ in the integral. Many
conjectures of this sort can be verified through the combinatorial models. On the
other hand, the really outstanding conjecture about the integral is that it
exists! At the present time there is no measure theory or generalization of
measure theory that supports it.   Here is a formal structure of great beauty. It
is also a structure whose consequences can be verified by a remarkable variety of
alternative means. \vspace{3mm}

In this section we will examine the formalism of Witten's approach via a generalization
of our sketch of ``integration without integration".  In order to do this we need to consider 
functions $f(A)$ of gauge connections $A$ and a notion of equivalence, $f \sim g,$ taking the form
$f-g =Dh$ where $D$ is a gauge functional derivative. Since these notions need defining, we first discuss them
in the context of the integrand of Witten's integral. Thus for a while, we shall speak of Witten's integral, but 
let it be known that this integral will soon be replaced by an equivalence class of functions just as happened in 
the last section!
\bigbreak

The formalism of the Witten integral implicates invariants of knots and links corresponding to
each classical Lie algebra.   In order to see this, we need to introduce the Wilson
loop.  The Wilson loop is an exponentiated version of integrating the gauge field
along a loop  $K$  in three space that we take to be an embedding (knot) or a
curve with transversal self-intersections.  For this discussion, the Wilson loop
will be denoted by the notation  $$W_{K}(A) = <K|A>$$ to denote the dependence on
the loop $K$ and the field $A$.   It is usually indicated by the symbolism  
$tr(Pe^{\oint_{K} A})$ .   Thus $$W_{K}(A) = <K|A>  = tr(Pe^{\oint_{K} A}).$$   
Here the $P$  denotes  path ordered integration - we are integrating and
exponentiating matrix valued functions, and so must keep track of the order of
the operations.  The  symbol  $tr$  denotes the trace of the resulting matrix.
This Wilson loop integration exists by normal means and will not be replaced by function classes. 
\vspace{3mm}

With the help of the Wilson loop functional on knots and links,  Witten  writes
down a functional integral for link invariants in a 3-manifold  $M$:

$$Z(M,K) = \int DAe^{(ik/4 \pi)S(M,A)} tr(Pe^{\oint_{K} A}) $$

$$= \int DAe^{(ik/4 \pi)S}<K|A>.$$

\noindent Here $S(M,A)$ is the Chern-Simons Lagrangian, as in the previous
discussion. We abbreviate  $S(M,A)$  as $S$ and write  $<K|A>$  for the Wilson
loop. Unless otherwise mentioned, the manifold  $M$  will be the
three-dimensional sphere  $S^{3}$ \vspace{3mm}

An analysis of the formalism of this functional integral reveals quite a bit
about its role in knot theory.   This analysis depends upon key facts relating
the curvature of the gauge field to both the Wilson loop and the Chern-Simons
Lagrangian. The idea for using the curvature in this way is due to Lee Smolin
\cite{Smolin} (See also \cite{Ramusino}). To this end, let us recall the local
coordinate structure of the gauge field  $A(x)$,  where $x$  is a point in
three-space.   We can write   $A(x)  =  A^{a}_{k}(x)T_{a}dx^{k}$  where  the
index  $a$ ranges from $1$ to $m$ with the Lie algebra basis $\{T_{1}, T_{2},
T_{3}, ..., T_{m}\}$.  The index $k$   goes from $1$  to $3$.     For each choice
of $a$  and  $k$,  $A^{a}_{k}(x)$   is a smooth function defined on three-space.
In  $A(x)$  we sum over the values of repeated indices.  The Lie algebra
generators $T_{a}$  are  matrices  corresponding to a given representation of the
Lie algebra of the gauge group $G.$   We assume some properties of these matrices
as follows: \vspace{3mm}

\noindent 1.  $[T_{a} , T_{b}] = i f^{abc}T_{c}$  where  $[x ,y] = xy - yx$ , and
$f^{abc}$ (the matrix of structure constants)  is totally antisymmetric.  There
is summation over repeated indices. \vspace{3mm}

\noindent 2.  $tr(T_{a}T_{b}) = \delta_{ab}/2$ where  $\delta_{ab}$ is the
Kronecker delta  ($\delta_{ab} = 1$ if $a=b$ and zero otherwise). \vspace{6mm}

We also assume some facts about curvature. (The reader may enjoy comparing with
the exposition in \cite{K and P}.  But note the difference of conventions on the
use of  $i$ in the Wilson loops and curvature definitions.)   The first fact  is
the relation of Wilson loops and curvature for small loops: \vspace{3mm}

\noindent {\bf Fact 1.} The result of evaluating a Wilson loop about a very small
planar circle around a point $x$ is proportional to the area enclosed by this
circle times the corresponding value of the curvature tensor of the gauge field
evaluated at $x$. The curvature tensor is  written 
$$F^{a}_{rs}(x)T_{a}dx^{r}dy^{s}.$$ It is the local coordinate expression of  $F
= dA +A \wedge A.$ \vspace{3mm}

\noindent {\bf Application of Fact 1.}  Consider a given Wilson line  $<K|S>$.
Ask how its value will change if it is deformed infinitesimally in the
neighborhood of a point $x$ on the line.  Approximate the change according to
Fact 1, and regard the point $x$ as the place of curvature evaluation.  Let 
$\delta<K|A>$  denote the change in the value of the line.    $\delta <K|A>$  is
given by the formula $$\delta <K|A> = dx^{r}dx^{s}F_{a}^{rs}(x)T_{a}<K|A>.$$ This
 is the first order approximation to the change in the Wilson line. \vspace{3mm}

In this formula it  is understood that the Lie algebra matrices  $T_{a}$  are to
be inserted into the Wilson line at the point $x$,  and that we are summing over
repeated indices. This means that  each  $T_{a}<K|A>$ is  a new Wilson line
obtained from the original  line $<K|A>$  by leaving the form of the loop
unchanged,  but inserting the matrix  $T_{a}$ into that loop at the point  $x$.  
In Figure 6 we have illustrated this mode of insertion of Lie algebra into the
Wilson loop. Here and in further illustrations in this section we use $W_{K}(A)$
to denote the Wilson loop. Note that in the diagrammatic version shown in Figure
6 we have let small triangles with legs indicate $dx^{i}.$ The legs correspond to
indices just as in our work in the last section with Lie algebras and chord
diagrams. The curvature tensor is indicated as a circle with three legs
corresponding to the indices of $F_{a}^{rs}.$ \vspace{3mm}

\noindent {\bf Notation.} In the diagrams in this section we have dropped mention
of the factor of $(1/ 4 \pi)$ that occurs in the integral. This convention saves
space in the figures. In these figures $L$ denotes the Chern--Simons Lagrangian.
\vspace{3mm}

\begin{figure}[htb]
     \begin{center}
     \begin{tabular}{c}
     \includegraphics[width=6cm]{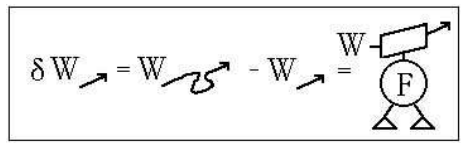}
     \end{tabular}
     \caption{\bf Lie Algebra and Curvature Tensor Insertion into the Wilson Loop}
     \label{Figure 6}
\end{center}
\end{figure}

\noindent {\bf Remark.}  In thinking about the Wilson line $<K|A> =
tr(Pe^{\oint_{K} A})$,  it is helpful to recall Euler's formula for the
exponential:

$$e^{x} = lim_{n \rightarrow \infty}(1+x/n)^{n}.$$

\noindent The  Wilson line  is  the limit, over partitions of the loop $K$,  of
products of the matrices  $(1 + A(x))$  where $x$ runs over the partition.  Thus
we can write symbolically,

$$<K|A> =  \prod_{x \in K}(1 +A(x))$$ $$=  \prod_{x \in K}(1 +
A^{a}_{k}(x)T_{a}dx^{k}).$$

\noindent It is understood that a product of matrices around a closed loop
connotes the trace of the product.  The ordering is forced by the one dimensional
nature of the loop.   Insertion of a given matrix into this product at a point on
the loop is then a well-defined concept.   If  $T$ is a given matrix then it is
understood that   $T<K|A>$  denotes the insertion of $T$ into some point of the
loop.  In the case above, it is understood from context in the formula that the
insertion is to be performed at the point $x$  indicated in the argument of the
curvature. \vspace{3mm}

\noindent {\bf Remark.}  The   previous remark implies the following formula for
the variation of the Wilson loop with respect to the gauge field:

$$\delta <K|A>/\delta (A^{a}_{k}(x))  =  dx^{k}T_{a}<K|A>.$$

\noindent Varying the Wilson loop with respect to the gauge field results in the
insertion of an infinitesimal Lie algebra element into the loop. Figure 7 gives a
diagrammatic form for this formula. In that Figure we use a capital $D$ with up
and down legs to denote the derivative $\delta /\delta (A^{a}_{k}(x)).$
Insertions in the Wilson line are indicated directly by matrix boxes placed in a
representative bit of line. \vspace{3mm}

\begin{figure}[htb]
     \begin{center}
     \begin{tabular}{c}
     \includegraphics[width=6cm]{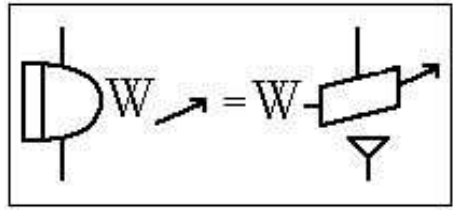}
     \end{tabular}
     \caption{\bf Differentiating the Wilson Line}
     \label{Figure 7}
\end{center}
\end{figure}

\noindent {\bf Proof.} $$\delta <K|A>/\delta (A^{a}_{k}(x))$$

$$= \delta \prod_{y \in K}(1 + A^{a}_{k}(y)T_{a}dy^{k})/\delta (A^{a}_{k}(x))$$

$$= \prod_{y<x \in K}(1 + A^{a}_{k}(y)T_{a}dy^{k}) [T_{a}dx^{k}] \prod_{y>x \in
K}(1 + A^{a}_{k}(y)T_{a}dy^{k})$$

$$= dx^{k}T_{a}<K|A>.$$ \vspace{3mm}

\noindent {\bf Fact 2.}  The variation of the Chern-Simons Lagrangian  $S$  with
respect to the gauge potential at a given point in three-space is related to the
values of the curvature tensor at that point by the following formula:

$$F^{a}_{rs}(x)  =  \epsilon_{rst} \delta S/\delta (A^{a}_{t}(x)).$$

\noindent Here  $\epsilon_{abc}$ is the epsilon symbol for three indices, i.e. it
is $+1$ for positive permutations of $123$ and $-1$ for negative permutations of
$123$ and zero if any two indices are repeated.  A diagrammatic for this formula
is shown in Figure 8. \vspace{3mm}

\begin{figure}[htb]
     \begin{center}
     \begin{tabular}{c}
     \includegraphics[width=6cm]{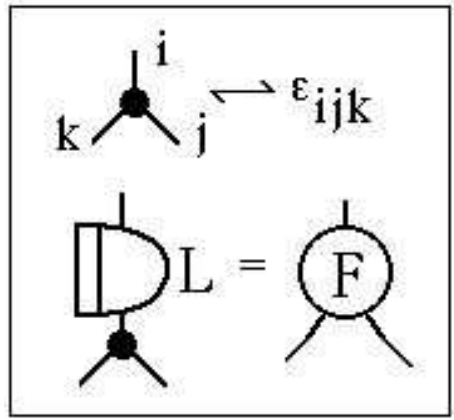}
     \end{tabular}
     \caption{\bf Variational Formula for Curvature}
     \label{Figure 8}
\end{center}
\end{figure}

\noindent {\bf The Functional Equivalence Relation.} With these facts at hand, we are prepared to define our equivalence
relation on functions of gauge fields. Given a function $F(A)$ of a gauge field $A$, we let
$DF$ denote any gauge functional derivative of $f(A).$ That is
$$DF = \delta F(A)/\delta (A^{a}_{k}(x)).$$ Note that 

$$D<K|A> = \delta <K|A>/\delta (A^{a}_{k}(x))  =  dx^{k}T_{a}<K|A>$$ with the insertion conventions as 
explained above. Then we say that functionals $F$ and $G$ are {\it integrally equivalent} ($F \sim G$) if there exists an $H$ 
such that $DH = F - G.$ We stipulate that all functionals in the discussion are rapidly vanishing at infinity,
where this is taken to mean that $F(A)$ goes to zero as $||A||$ goes to infinity, and the same is true for all 
functional derivatives of $F.$ Here the norm $$||A|| = \Sigma_{i,a} \int_{R^{3}}(A^{a}_{i})^{2}dvol$$ where $dvol$ is
the volume form on $R^3$ and it is assumed that all gauge fields have finite norm in this sense. 
\bigbreak

\noindent We then define the integral
$$Z(M,K) = \int DAe^{(ik/4 \pi)S(M,A)} tr(Pe^{\oint_{K} A}) = \int DAe^{(ik/4 \pi)S}<K|A>$$ to be the equivlance 
class of the functional $$e^{(ik/4 \pi)S(M,A)} tr(Pe^{\oint_{K} A}).$$ We invite the reader to make this
interpretation throughout the derivations that follow. It will then be apparent that much of what is usually taken for
formal heuristics about the funtional integral is actually a series of structural remarks about these equivalence
classes. Of course, one needs to know that the equivalence classes are non-trivial to make a complete story.
An existent integral would supply that key ingredient.  It its absence, we can examine that structure that can be
articulated at the level of the equivalence classes.
\bigbreak

We are prepared to determine how the Witten integral
behaves under a small deformation of the loop $K.$ \vspace{3mm}

\noindent {\bf Theorem.} 1. Let   $Z(K) = Z(S^{3},K)$  and let  $\delta Z(K)$ 
denote the change of $Z(K)$ under an infinitesimal change in the loop  K.   Then

$$ \delta Z(K) = (4 \pi i/k) \int dA e^{(ik/4\pi)S}[Vol] T_{a} T_{a} <K|A>$$

\noindent where $Vol = \epsilon_{rst} dx^{r} dx^{s} dx^{t}.$

The sum is taken over repeated indices, and the insertion is taken of the
matrices $T_{a}T_{a}$  at the chosen point  $x$  on the loop $K$ that is regarded
as the ÒcenterÓ of the deformation.  The volume element $Vol =
\epsilon_{rst}dx_{r}dx_{s}dx_{t}$ is taken with regard to the infinitesimal
directions of the loop deformation from this point on the original loop.
\vspace{3mm}

\noindent 2. The same formula applies, with a different interpretation,  to the
case where  $x$  is  a double point of transversal self intersection of a loop K,
 and the deformation consists in shifting one of the crossing segments
perpendicularly to the plane of intersection  so that the self-intersection point
disappears.  In this  case,  one  $T_{a}$  is inserted into each of the
transversal crossing segments so that  $T_{a}T_{a}<K|A>$ denotes a Wilson loop
with a self intersection  at  $x$   and insertions of $T_{a}$  at  $x +
\epsilon_{1}$ and  $x + \epsilon_{2}$  where $\epsilon_{1}$ and $\epsilon_{2}$
denote small displacements along the two arcs of $K$ that intersect at $x.$  In
this case, the volume form is nonzero, with two  directions coming from the plane
of movement of one arc, and the perpendicular direction is the direction of the
other arc. \vspace{3mm}

\noindent {\bf Proof.}

$$\delta Z(K)  =  \int DA e^{(ik/4 \pi)S} \delta <K|A>$$

$$= \int DA e^{(ik/4 \pi)S} dx^{r}dy^{s} F^{a}_{rs}(x) T_{a}<K|A>$$

$$=  \int DA e^{(ik/4 \pi)S} dx^{r}dy^{s} \epsilon_{rst} (\delta S/\delta
(A^{a}_{t}(x)))  T_{a}<K|A>$$

$$= (-4 \pi i/k) \int DA  (\delta e^{(ik/4 \pi)S}/\delta (A^{a}_{t}(x)))
\epsilon_{rst} dx^{r}dy^{s}T_{a}<K|A>$$

$$= (4 \pi i/k) \int DA  e^{(ik/4 \pi)S} \epsilon_{rst} dx^{r}dy^{s} (\delta
T_{a}<K|A>/\delta (A^{a}_{t}(x)))$$

(integration by parts and the boundary terms vanish)

$$=  (4 \pi i/k) \int DA  e^{(ik/4 \pi)S}[Vol] T_{a}T_{a}<K|A>.$$

This completes the formalism of the proof.  In the case of part 2., a change of
interpretation occurs at the point in the argument when the Wilson line is
differentiated.  Differentiating a self-intersecting Wilson line at a point of
self intersection is equivalent to differentiating the corresponding product of
matrices with respect to a variable that occurs at two points in the product
(corresponding to the two places where the loop passes through the point).  One
of these derivatives gives rise to a term with volume form equal to zero, the
other term is the one that is described in part 2.  This completes the proof of
the Theorem.  // \vspace{3mm}

\noindent The formalism of this proof is illustrated in Figure 9. \vspace{3mm}

\begin{figure}[htb]
     \begin{center}
     \begin{tabular}{c}
     \includegraphics[width=6cm]{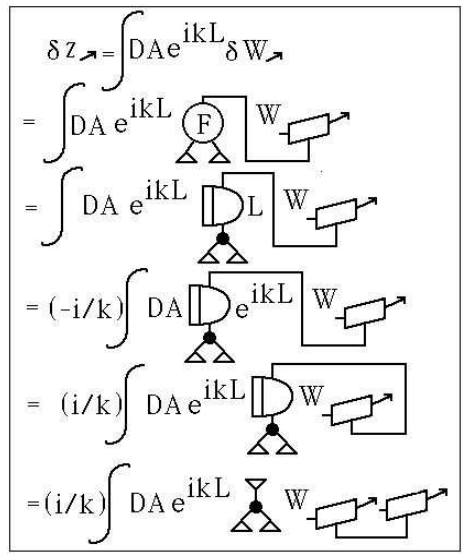}
     \end{tabular}
     \caption{\bf Varying the Functional Integral by Varying the Line}
     \label{Figure 9}
\end{center}
\end{figure}

In the case of switching a crossing the key point is to write the crossing switch
as a composition of first moving a segment to obtain a transversal intersection
of the diagram with itself, and then to continue the motion to complete the
switch.  One then analyzes separately the case where $x$  is  a double point of
transversal self intersection of a loop $K,$  and the deformation consists in
shifting one of the crossing segments perpendicularly to the plane of
intersection  so that the self-intersection point disappears.  In this  case, 
one  $T_{a}$  is inserted into each of the transversal crossing segments so that 
$T^{a}T^{a}<K|A>$ denotes a Wilson loop with a self intersection  at  $x$   and
insertions of $T^{a}$  at $x + \epsilon_{1}$ and  $x + \epsilon_{2}$   as in part
$2.$ of the Theorem above. The first insertion is in the moving line, due to
curvature. The second insertion is the consequence of differentiating the
self-touching Wilson line. Since this line can be regarded as a product, the
differentiation occurs twice at the point of intersection, and it is the second
direction that produces the non-vanishing volume form. \vspace{3mm}

Up to the choice of our conventions for constants, the switching formula is, as
shown below (See Figure 10).

$$Z(K_{+}) -  Z(K_{-}) =  (4 \pi i/k)\int DA  e^{(ik/4\pi)S}
T_{a}T_{a}<K_{**}|A>$$ $$= (4 \pi i/k) Z(T^{a}T^{a}K_{**}),$$

\noindent where $K_{**}$ denotes the result of replacing the crossing by a
self-touching crossing. We distinguish this from adding a graphical node at this
crossing by using the double star notation. \vspace{3mm}

\begin{figure}[htb]
     \begin{center}
     \begin{tabular}{c}
     \includegraphics[width=6cm]{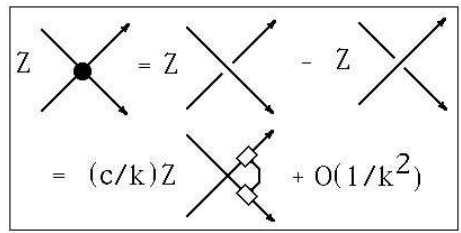}
     \end{tabular}
     \caption{\bf The Difference Formula}
     \label{Figure 10}
\end{center}
\end{figure}

A key point is to notice that the Lie algebra insertion for this difference is
exactly what is done (in chord diagrams) to make the weight systems for Vassiliev
invariants (without the framing compensation).  In order to extend the Heuristic at this point we need to 
assume the analog of a perturbative expansion for the integral. That is, we assume that that there are
invariants of regular isotopy of $K$,  $Z_{n}(K)$ and that 
$$ e^{(ik/4 \pi)S(A)}<K|A>~ \sim~ \Sigma_{n=0}^{\infty} k^{-n} Z_{n}(K).$$ Note that since we have shown that the
equivalence class of $$ e^{(ik/4 \pi)S(A)}<K|A> $$ is a regular isotopy invariant, it is not at all implausible to assume
that there is a power series representative of this functional whose coefficients are numerical regular isotopy
invariants. It is this assumption that allows one to make contact with numerical evaluations. The assumption of this
power series  representation corresponds to the formal
perturbative expansion of the Witten integral. One obtains Vassiliev invariants as
coefficients of the powers of ($1/k^{n}$). Thus the formalism of the Witten
functional integral takes one directly to these weight systems in the case of the
classical Lie algebras. In this way the functional integral  is central to the
structure of the Vassiliev invariants. \vspace{3mm}

\section{Perturbative Expansion}
Letting $M^3$ be a three-manifold and $K$ a knot or link in $M^{3},$
we write $$\psi(A) = e^{ik L(A)}W_{K}(A),$$ and replace $A$ by $A/\sqrt{k}$  then we can write
$$\hat{\psi}(A) = e^{\frac{i}{4 \pi} \int_{M^3} tr(A \wedge dA)} e^{ \frac{i}{6\pi \sqrt{k}} \int_{M^3} tr(A \wedge A \wedge A) }W_{K}(A/\sqrt{k}).$$
It is the equivalence class of this functional of gauge fields that contains much topological information about knots and links in the three-manifold $M^{3}.$
We can expand this functional by taking the explicit formula for the Wilson loop:
$$W_{K}(A/\sqrt{k}) =   tr(\prod_{{x}\in K}(1 + A({x})/\sqrt{k}).$$

$$
\hat{\psi}(A) =  e^{\frac {i}{4 \pi} \int_{M^3} tr(A \wedge dA)}
 e^{\frac{i}{6\pi}\int_{M^3} tr(\frac{1}{\sqrt k} A\wedge A\wedge A)}
 tr(\prod_{{x}\in K}(1 + A({x})/\sqrt k)).
 $$
 
 $$tr\Bigl(\prod_{{x}\in K}(1 + \frac1{\sqrt k}A({x}))\Bigr)
 =tr\Bigl(1 + \frac{1}{\sqrt k}\int\limits_K A +
  \frac{1}{k} \int\limits_{K_1 < K_2}A({x}_1)A({x}_2)+\ldots\Bigr)
$$ 
where
$$
 \int\limits_{K_1<\ldots<K_n} A({x}_1)A({x}_2)\ldots A({x}_n)=
 \int\limits_{\underbrace{\scriptstyle K\times\ldots\times K}_{n}=K^n}A({x}_1)A({x}_2)\ldots
 A({x}_n)
 $$
 $$\overrightarrow{{x}} = ({x}_1,{x}_2,\ldots ,{x}_n)\in K^n\ \hbox{with}\
 {x}_1<{x}_2<\ldots <{x}_n.$$
This is an iterated integrals expression for the Wilson loop.\\

Our functional is transformed into a perturbative series in powers of $1/k .$ The equivalence class of each term in the series (when $M^3$ is the three-sphere $S^3$) is formally a Vassiliev invariant as we have described in the 
previous section. A more intense look at the structure of these functionals can be accomplished by gauge-fixing as we show in the last section.\\

\section{The Loop Transform and Loop Quantum Gravity}

Suppose that $\psi (A)$ is a (complex-valued) function defined on gauge fields. Then we define formally the {\em loop transform}
$\widehat{\psi}(K)$, a function on embedded loops in three dimensional space, by the formula

$$\widehat{\psi}(K) = \int \psi(A) W_{K}(A).$$

\noindent note that we could also write
$$\widehat{\psi}(K) = \psi(A) W_{K}(A)$$
where it is understood that the right-hand side of the equation represents its integral equivalence class. Then we can look at it as a function of the loop $K$ {\it and} as a function of the 
gauge field $A$. This changes one's point of view about the loop transform. We are really examining a hybrid function of both a possibly knotted loop $K$ and a gauge field $A.$ The important 
structure is the relationships that ensue in the integral equivalence class between varying $A$ and varying $K.$ Nevertheless, we shall continue to use integral signs to remind the reader that we are working
with the integral equivalence classes of these functionals.\\ 

\noindent
If $\Delta$ is a differential operator defined on $\psi(A),$ then we can use this integral transform to shift the effect of $\Delta$ to an operator on loops via integration by parts:

$$\widehat{ \Delta \psi }(K) = \int \Delta \psi(A) W_{K}(A)$$

$$ = - \int   \psi(A) \Delta W_{K}(A).$$

\noindent
When $\Delta$ is applied to the Wilson loop the result can be an understandable geometric or topological operation.  In Figures 11, 12 and 13 we illustrate this situation with diagrammatically defined operators $G$ and $H.$ 
\vspace{3mm}

 \begin{figure}[htb]
     \begin{center}
     \begin{tabular}{c}
     \includegraphics[width=6cm]{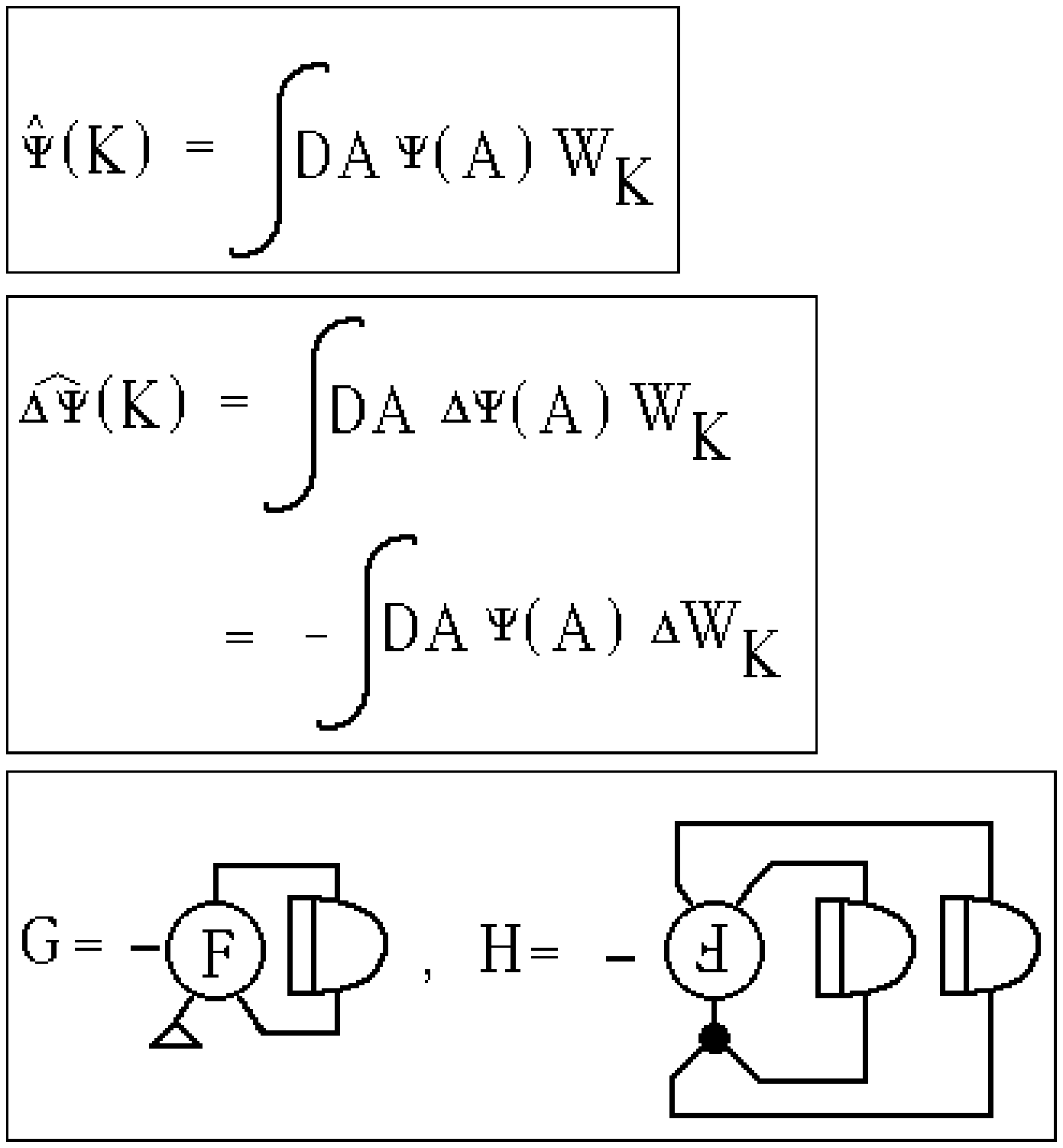}
     \end{tabular}
     \caption{\bf The Loop Transform and Operators $G$ and $H$}
     \label{Figure 11}
\end{center}
\end{figure}

 \begin{figure}[htb]
     \begin{center}
     \begin{tabular}{c}
     \includegraphics[width=6cm]{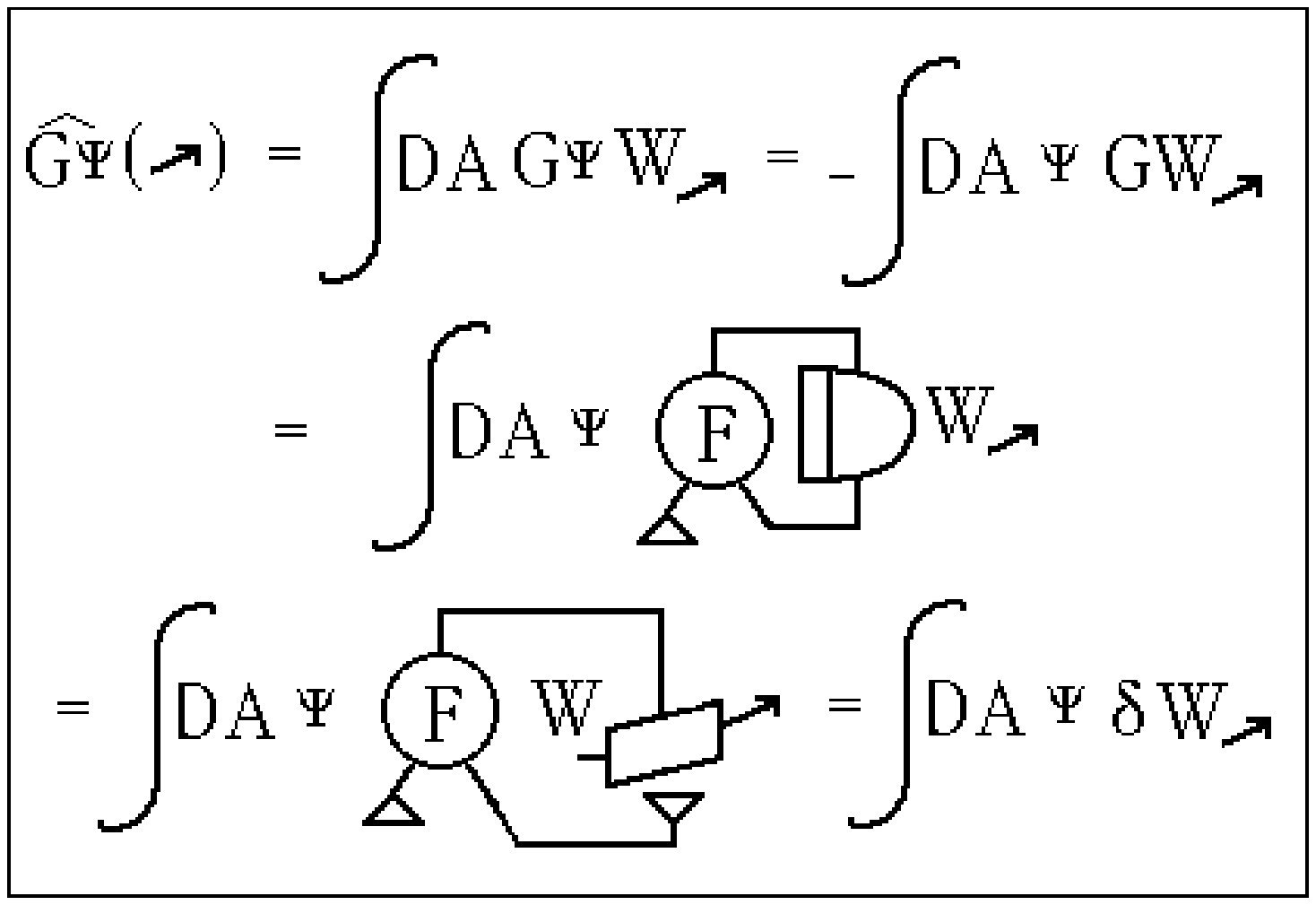}
     \end{tabular}
     \caption{\bf The Diffeomorphism Constraint}
     \label{Figure 12}
\end{center}
\end{figure}

 \begin{figure}[htb]
     \begin{center}
     \begin{tabular}{c}
     \includegraphics[width=6cm]{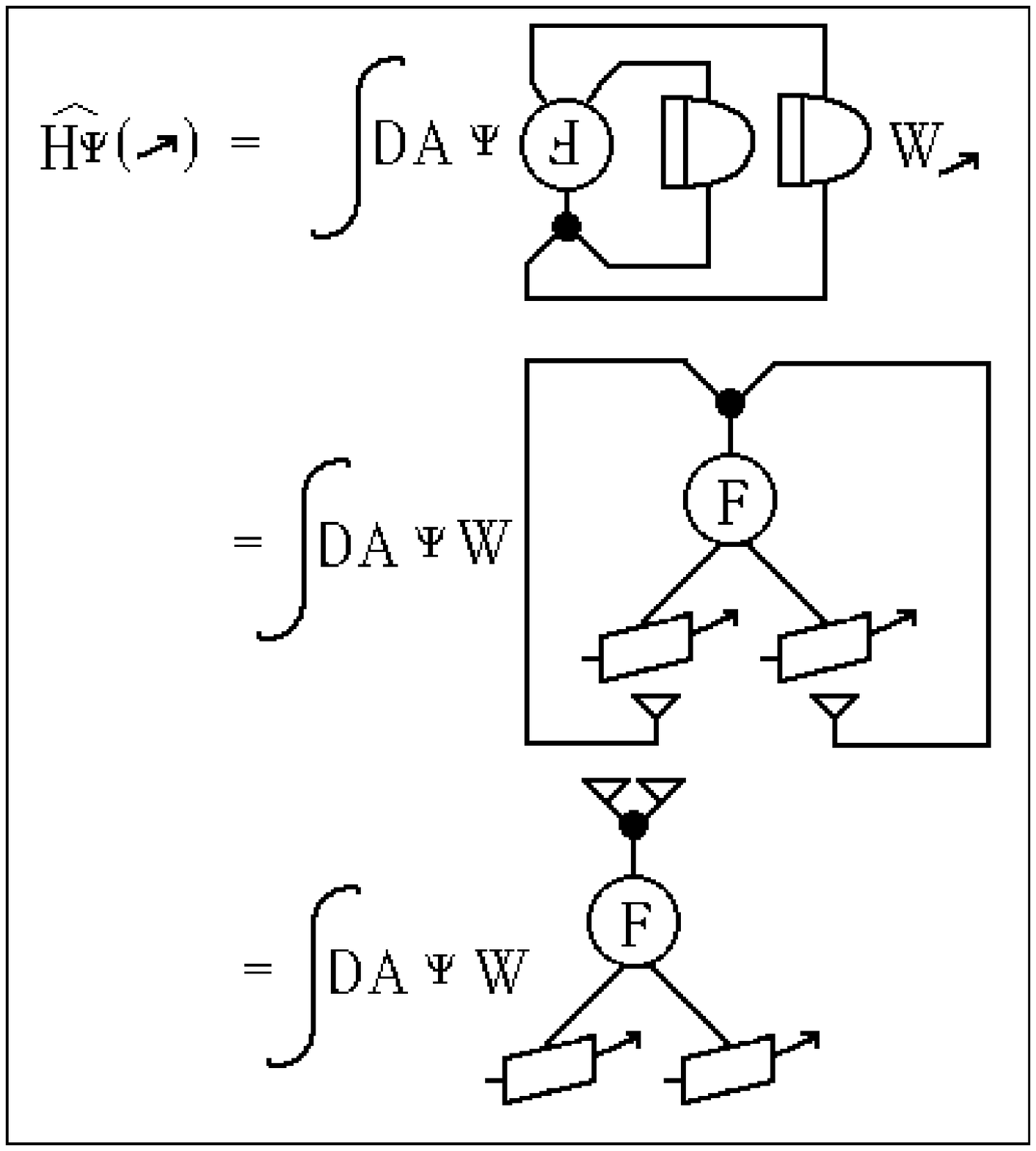}
     \end{tabular}
     \caption{\bf The Hamiltonian Constraint}
     \label{Figure 13}
\end{center}
\end{figure}

\noindent
We see from Figure 12 that 

$$\widehat{ G \psi }(K) =  \delta \widehat{ \psi }(K)$$

\noindent
where this variation refers to the effect of  varying $K$ by a small loop. As we saw in this section, this means that if  $\widehat{ \psi }(K)$  is a topological invariant of knots and links, then 
$\widehat{ G \psi }(K) =0$ for all embedded loops $K.$  This condition is a transform analogue of the equation $G \psi(A) =0.$
This equation is the differential analogue of an invariant of knots and links. It may happen that $\delta \widehat{ \psi }(K)$ is not strictly zero, as in the case of our framed knot invariants.
For example with $$\psi(A) =  e^{(ik/4\pi) \int tr(A \wedge dA + (2/3)A \wedge A \wedge A)}$$
we conclude that $\widehat{ G \psi }(K)$ is zero for flat deformations (in the sense of this section) of the loop $K,$ but can be non-zero in the presence of a twist or curl. In this sense the loop transform provides a subtle variation on the strict condition $G \psi(A) =0.$ This Chern-Simons functional $\psi(A)$ can be seen to be a state of loop quantum gravity.
\vspace{3mm}

In \cite{ASR} and earlier publications by these authors, the loop transform is used to study a reformulation and quantization of Einstein gravity.  The differential geometric gravity theory is reformulated in terms of a background gauge connection and in the quantization, the Hilbert space consists in functions $\psi(A)$ that are required to satisfy the constraints 
$$G \psi =0$$
and 
$$H \psi =0.$$

\noindent
where $H$ is the operator shown in Figure 13. Thus we see that 
$\widehat{G}(K)$ can be partially zero in the sense of producing a framed knot invariant, and (from Figure 13 and the antisymmetry of the epsilon) that $\widehat{H}(K)$ is zero for non-self-intersecting loops.  This means that the loop transforms of $G$ and $H$ can be used to investigate a subtle variation of the original scheme for the quantization of gravity.  The appearance of the Chern-Simons state
$$\psi(A) =  e^{(ik/4\pi) \int tr(A \wedge dA + (2/3)A \wedge A \wedge A)}$$
is quiite remarkable in this theory, where it is commonly referred to as the {\it Kodama State.} See \cite{Smolin1,Maguelijo,Wieland,Randono,Kodama,Soo,Pullin,WittKodama}  for a number of references about this 
state, up to the present day. Many ways of weaving this relationship of knot theory and quantum gravity have been devised, from examining directly the Kodama state and its relationship with DeSitter space, to
the evolution of spin networks and spin foams to handle the fundamental topological conditions in the theory. 
\vspace{3mm}

\section{Wilson Lines, Axial Gauge and the Kontsevich Integrals}

In this section we follow the gauge fixing method used by Fr\"ohlich and King 
\cite{Frohlich and King}.  Their paper was written before the advent of Vassiliev 
invariants, but contains, as we shall see, nearly the whole story about the Kontsevich 
integral.  A similar approach to ours can be found in \cite{LP}. In our case we have simplified the determination of the inverse operator for this formalism and we have given a few more details about the calculation of the correlation functions than is customary in physics literature.  I hope that this approach makes this subject more accessible to mathematicians. A heuristic argument of this kind contains a great deal of valuable mathematics. It is clear that these matters will eventually be given a fully rigorous treatment. In fact, in the present case there is a rigorous treatment, due to Albevario and Sen-Gupta \cite{AS} of the functional integral {\em after} the light-cone  gauge has been imposed.
\vspace{3mm}

\noindent
Let $(x^{0}, x^{1}, x^{2})$ denote a point in three dimensional space.
Change to light-cone coordinates
$$x^{+} = x^{1} + x^{2}$$ and
$$x^{-} = x^{1} - x^{2}.$$

\noindent
Let $t$ denote $x^{0}.$
\vspace{3mm}

\noindent
Then the gauge connection can be written in the form
$$A(x) = A_{+}(x)dx^{+} + A_{-}(x)dx^{-} + A_{0}(x)dt.$$
\vspace{3mm}

\noindent
Let $CS(A)$ denote the Chern-Simons integral (over the three dimensional sphere)

$$CS(A) = (1/4\pi)\int tr(A \wedge dA + (2/3) A \wedge A \wedge A).$$

\noindent
We define {\em axial gauge} to be the condition that $A_{-} = 0.$
We shall now work with the functional integral of the previous section under the axial 
gauge restriction. In axial gauge we have that 
$A \wedge A \wedge A = 0$ and so 

$$CS(A) = (1/4\pi)\int tr(A \wedge dA).$$

\noindent
Letting $\partial_{\pm}$ denote partial differentiation with respect to $x^{\pm}$, we get 
the following formula in axial gauge
$$A \wedge dA = (A_{+} \partial_{-} A_{0} - A_{0} \partial_{-}A_{+})dx^{+} \wedge 
dx^{-} \wedge dt.$$

\noindent
Thus, after integration by parts, we obtain the following formula for the Chern-Simons 
integral:

 $$CS(A) = (1/2 \pi) \int tr(A_{+} \partial_{-} A_{0}) dx^{+} \wedge dx^{-} \wedge 
dt.$$

\noindent
Letting $\partial_{i}$ denote the partial derivative with respect to $x_{i}$, we have that 
$$\partial_{+} \partial_{-} = \partial_{1}^{2} - \partial_{2}^{2}.$$  If we replace 
$x^{2}$ with $ix^{2}$ where 
$i^{2} = -1$, then $\partial_{+} \partial_{-}$  is replaced by
$$\partial_{1}^{2} + \partial_{2}^{2} = \nabla^{2}.$$ 
We now make this replacement so that the analysis can be expressed over the complex 
numbers.
\vspace{3mm}

\noindent
Letting $$z = x^{1} + ix^{2},$$ it is well known that 
$$\nabla^{2} ln(z) = 2 \pi \delta(z)$$
where $\delta(z)$ denotes the Dirac delta function and $ln(z)$ is the natural logarithm of 
$z.$  Thus we can write
$$(\partial_{+} \partial_{-})^{-1} = (1/2 \pi)ln(z).$$
Note that $\partial_{+} = \partial_{z} = \partial /\partial z$  after the replacement of 
$x^{2}$ by $ix^{2}.$   As a result we have that 
$$(\partial_{-})^{-1} = \partial_{+} (\partial_{+} \partial_{-})^{-1} =
\partial_{+} (1/2 \pi)ln(z) = 1/2 \pi z.$$

\noindent
Now that we know the inverse of the operator $\partial_{-}$ we are in a position to treat 
the Chern-Simons integral as a quadratic form in the pattern 
$$ (-1/2)<A, LA> = - iCS(A)$$  
where the operator 
$$L =  \partial_{-}.$$
Since we know $L^{-1}$, we can express the functional integral as a Gaussian integral:
\vspace{3mm}

\noindent
We replace 

$$Z(K) = \int DAe^{ikCS(A)} tr(Pe^{\oint_{K} A})$$ 

\noindent
by

$$Z(K) = \int DAe^{iCS(A)} tr(Pe^{\oint_{K} A/\sqrt k})$$ 

\noindent
by sending $A$ to $(1/ \sqrt k)A$. We then replace this version by

$$Z(K) = \int DAe^{(-1/2)<A, LA>} tr(Pe^{\oint_{K} A/\sqrt k}).$$

\noindent
In this last formulation we can use our knowledge of $L^{-1}$ to determine the the 
correlation functions and express $Z(K)$ perturbatively in powers of $(1/ \sqrt k).$
\vspace{3mm}

\noindent
{\bf Proposition.}
Letting 
$$<\phi(A)> = \int DA e^{(-1/2)<A, LA>}\phi(A) / \int DA e^{(-1/2)<A, LA>}$$  for 
any functional $\phi(A)$,
we find that
$$<A_{+}^{a}(z,t)A_{+}^{b}(w,s)> = 0,$$
$$<A_{0}^{a}(z,t)A_{0}^{b}(w,s)> = 0,$$
$$<A_{+}^{a}(z,t)A_{0}^{b}(w,s)> = \kappa \delta^{ab} \delta(t-s)/(z-w)$$  where 
$\kappa$ is a constant.
\vspace{3mm}

\noindent
{\bf Proof Sketch.}
Let's recall how these correlation functions are obtained.
The basic formalism for the Gaussian integration is in the pattern

$$<A(z)A(w)> = \int DA e^{(-1/2)<A, LA>} A(z)A(w) / \int DA e^{(-1/2)<A, LA>}$$

$$ =  ((\partial / \partial J(z)) (\partial / \partial J(w)) |_{J=0})
 e^{(1/2)<J, L^{-1}J>}$$

\noindent
Letting $G*J(z) = \int dw G(z-w)J(w)$, we have that when

$$LG(z) = \delta(z)$$  

\noindent
($\delta(z)$ is a Dirac delta function of $z$.) then

$$LG*J(z) = \int dw LG(z-w)J(w) = \int dw \delta(z-w) J(w) = J(z)$$

\noindent
Thus $G*J(z)$ can be identified with $L^{-1}J(z)$.
\vspace{3mm}

\noindent
In our case 
$$G(z) = 1/ 2 \pi z$$  

\noindent
and 
 
$$L^{-1}J(z) = G*J(z) = \int dw J(w)/(z-w).$$

\noindent
Thus 

$$<J(z),L^{-1}J(z)> = <J(z), G*J(z)> = (1/ 2\pi) \int tr(J(z) (\int dw J(w)/(z-w)) dz$$

$$ = (1/ 2\pi) \int \int dz dw tr(J(z)J(w))/(z-w).$$

\noindent
The results on the correlation functions then follow directly from differentiating this 
expression.  Note that the Kronecker delta on Lie algebra indices is a result of the corresponding Kronecker delta in the trace formula
$tr(T_{a}T_{b}) = \delta_{ab}/2$ for products of Lie algebra generators. The Kronecker delta for the $x^{0} = t, s$ coordinates is a consequence of the evaluation at $J$ equal to zero.//
\vspace{3mm}

We are now prepared to give an explicit form to the perturbative expansion for 

$$<K>= Z(K)/\int DAe^{(-1/2)<A, LA>}$$

$$= \int DAe^{(-1/2)<A, LA>} tr(Pe^{\oint_{K} A/\sqrt k})/ \int DAe^{(-1/2)<A, LA>}$$

$$ = \int DAe^{(-1/2)<A, LA>} tr(\prod_{x \in K} (1 +  (A/\sqrt k)))/\int DAe^{(-
1/2)<A, LA>}$$

$$ = \sum_{n} (1/k^{n/2}) \oint_{K_{1} < ... < K_{n}} <A(x_{1}) ... A(x_{n})>.$$

\noindent
The latter summation can be rewritten (Wick expansion) into a sum over products of pair 
correlations, and we have already worked out the values of these. In the formula above we 
have written  $K_{1} < ... < K_{n}$ to denote the integration over variables $x_{1} , ... 
x_{n}$ on $K$ so that $x_{1} <  ... < x_{n}$ in the ordering induced on the loop $K$ by 
choosing a basepoint on the loop.  After the Wick expansion, we get

$$<K> = \sum_{m} (1/k^{m}) \oint_{K_{1} < ... < K_{n}} 
\sum_{P= \{x_{i} < x'_{i}| i = 1, ... m\}}
\prod_{i}<A(x_{i})A(x'_{i})>.$$

\noindent
Now we know that 

$$<A(x_{i})A(x'_{i})> = 
<A^{a}_{k}(x_{i})A^{b}_{l}(x'_{i})>T_{a}T_{b}dx^{k}dx^{l}.$$

\noindent
Rewriting this in the complexified axial gauge coordinates, the only contribution is 

$$<A_{+}^{a}(z,t)A_{0}^{b}(s,w)> = \kappa \delta^{ab} \delta(t-s)/(z-w).$$

\noindent
Thus

$$<A(x_{i})A(x'_{i})>$$ 
$$=<A^{a}_{+}(x_{i})A^{a}_{0}(x'_{i})>T_{a}T_{a}dx^{+} \wedge dt + 
<A^{a}_{0}(x_{i})A^{a}_{+}(x'_{i})>T_{a}T_{a}dx^{+} \wedge dt$$
$$= (dz-dz')/(z-z') [i/i']$$
where $[i/i']$ denotes the insertion of the Lie algebra elements
$T_{a}T_{a}$ into the Wilson loop.
\vspace{3mm}

\noindent
As a result, for each partition of the loop and choice of pairings 
$P= \{x_{i} < x'_{i}| i = 1, ... m\}$ we get an evaluation $D_{P}$ of the trace of these 
insertions into the loop. This is the value of the corresponding chord diagram in the weight 
systems for Vassiliev invariants. These chord diagram evaluations then figure in our 
formula as shown below:

$$<K> = \sum_{m} (1/k^{m}) \sum_{P} D_{P} \oint_{K_{1} < ... < K_{n}}
\bigwedge_{i=1}^{m}(dz_{i} - dz'_{i})/((z_{i} - z'_{i})$$

\noindent
This is a Wilson loop ordering version of the Kontsevich integral. To see the usual form of 
the integral appear, we change from the time variable (parametrization) associated with the 
loop itself to  time variables associated with a specific global direction of time in three 
dimensional space that is perpendicular to the complex plane defined by the axial gauge 
coordinates. It is easy to see that this results in one change of sign for each segment of the 
knot diagram supporting a pair correlation where the segment is oriented (Wilson loop 
parameter)  downward with respect to the global time direction. This results in the rewrite 
of our formula to

$$<K> = \sum_{m} (1/k^{m}) \sum_{P} (-1)^{|P \downarrow |} D_{P} \int_{t_{1} < 
... < t_{n}}
\bigwedge_{i=1}^{m}(dz_{i} - dz'_{i})/((z_{i} - z'_{i})$$

\noindent
where $|P \downarrow |$ denotes the number of points $(z_{i},t_{i})$ or $(z'_{i},t_{i})$  
in the pairings where the knot diagram is  oriented downward with respect to global time. 
The integration around the Wilson loop has been replaced by integration in the vertical time 
direction and is so indicated by the replacement of $\{ K_{1} < ... < K_{n} \}$ with $\{ 
t_{1} < ... < t_{n} \}$
\vspace{3mm}

\noindent
The coefficients of $1/k^{m}$  in this expansion are exactly the Kontsevich integrals for 
the weight systems $D_{P}$. See Figure 14.
\vspace{3mm}

 \begin{figure}[htb]
     \begin{center}
     \begin{tabular}{c}
     \includegraphics[width=6cm]{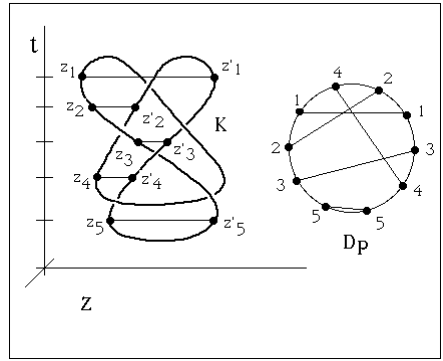}
     \end{tabular}
     \caption{\bf Applying the Kontsevich Integral}
     \label{Figure 14}
\end{center}
\end{figure}

\noindent
It was Kontsevich's insight to see (by different means) that 
these integrals could be used to construct Vassiliev invariants from arbitrary weight 
systems satisfying the four-term relations.  Here we have seen how these integrals arise 
naturally in the axial gauge fixing of the Witten functional integral. 
\vspace{3mm}

\noindent
{\bf Remark.}  The reader will note that we have not made a discussion of the role of the maxima and minima of the space curve of the knot with respect to the height direction ($t$).  In fact one has to take these maxima and minima very carefully into account and to divide by the corresponding evaluated loop pattern (with these maxima and minima) to make the Kontsevich integral well-defined and actually invariant under ambient  isotopy (with appropriate framing correction as well). The corresponding difficulty appears here in the fact that because of the gauge choice the Wilson lines are actually only defined in the complement of the maxima and minima and one needs to analyse a limiting procedure to take care of the inclusion of these points in the Wilson line.\\


\begin{thebibliography}{10}
\bibitem{AS}
S. Albevario and A. Sen-Gupta,
A Mathematical Construction of the Non-Abelian Chern-Simons Functional Integral,
{\em Commun. Math. Phys.}, Vol. 186 (1997), pp. 563-579.

\bibitem{ASR}
Ashtekar,Abhay, Rovelli, Carlo and Smolin,Lee [1992], "Weaving a Classical Geometry with Quantum Threads", {\em Phys. Rev. Lett.}, vol. 69, p. 237. 


\bibitem{Altschuler-Friedel}
Daniel Altschuler and Laurent Freidel,
Vassiliev knot invariants and Chern-Simons perturbation theory to all orders, {\em 
Commun. Math. Phys.} 187 (1997), 261-287.

\bibitem{Atiyah} M.F. Atiyah, {\em The Geometry and Physics of Knots}, Cambridge
University Press, 1990.

\bibitem{Bar-Natan} D. Bar-Natan,  On the Vassiliev knot invariants, {\em
Topology} {\bf 34} (1995), 423-472.

\bibitem{Bar-Natan-Thesis} Dror Bar-Natan, {\em Perturbative Aspects of the
Chern-Simons Topological Quantum field Theory}, Ph. D. Thesis, Princeton
University, June 1991.

\bibitem{Birman and Lin} J. Birman and X.S.Lin,  Knot polynomials and Vassiliev
invariants, {\em Invent. Math.} {\bf 111} No. 2 (1993), 225-270.

\bibitem{CDM} C. Dewitt-Morette, P. Cartier and A. Folacci, {\em Functional Integration - Basics and 
Applications}, NATO ASI Series, Series B: Physics Vol. 361 (1997).


\bibitem{Frohlich and King} J. Fr\"ohlich and C. King, The Chern-Simons Theory
and Knot Polynomials, {\em Commun. Math. Phys.} {\bf 126} (1989), 167-199.


\bibitem{HK} H. Kleinert, {\em Path Integrals in Quantum Mechanics, Statistics and Polymer Physics},
2nd edition, World Scientific, Singapore (1995).

\bibitem{Kauffman-Graph} L.H.Kauffman, New invariants in the theory of knots,
{\em Amer. Math. Monthly}, Vol.95,No.3,March 1988. pp 195-242.


\bibitem{Kauffman-Vogel} L.H.Kauffman and P.Vogel, Link polynomials and a
graphical calculus, {\em Journal of Knot Theory and Its Ramifications}, Vol. 1,
No. 1,March 1992, pp. 59- 104.

\bibitem{K and P} L.H.Kauffman, {\em Knots and Physics}, World Scientific 
Pub.,1991 and 1993

\bibitem{Func} L. H. Kauffman, Functional Integration and the theory of knots,
{J. Math. Physics}, Vol. 36 (5), May 1995, pp. 2402 - 2429.

\bibitem{WittKont} L. H. Kauffman, Witten's Integral and the Kontsevich
Integrals, in {\em Particles, Fields, and Gravitation}, Proceedings of the Lodz,
Poland (April 1998) Conference on Mathematical Physics edited by Jakub
Remblienski, AIP Conference Proceedings 453 (1998), pp. 368 -381.

\bibitem{Heuristics} L. H. Kauffman Knot Theory and the heuristics of functional integration,
{\em Physica A} 281 (2000), 173-200.

\bibitem{Kleinert} H. Kleinert, {\em Grand Treatise on Functional Integration}, World Scientific
 Pub. Co. (1999).

\bibitem{LP} J. M. F. Labastida and E. P$\acute{e}$rez, Kontsevich Integral for
Vassiliev Invariants from Chern-Simons Perturbation Theory in the Light-Cone
Gauge, {\em J. Math. Phys.}, Vol. 39 (1998), pp. 5183-5198.


\bibitem{Ramusino} P. Cotta-Ramusino,E.Guadagnini,M.Martellini,M.Mintchev,
Quantum field theory and link invariants, {\em Nucl. Phys. B} {\bf 330}, Nos. 2-3
(1990), pp. 557-574

\bibitem{Smolin} Lee Smolin, Link polynomials and critical points of the
Chern-Simons path integrals, {\em Mod. Phys. Lett. A}, Vol. 4,No. 12, 1989, pp.
1091-1112.

\bibitem{Smolin1}
Lee Smolin. Quantum gravity with a positive cosmological constant.
hep-th/0209079

\bibitem{Maguelijo}
Joao Magueijo, Laura Bethke. New ground state for quantum gravity.
arXiv:1207.0637

\bibitem{Wieland}
Wolfgang Wieland. Complex Ashtekar variables, the Kodama state and spinfoam gravity.
arXiv:1105.2330 

\bibitem{Randono}
Andrew Randono. In Search of Quantum de Sitter Space: Generalizing the Kodama State.  arXiv:0709.2905

\bibitem{Kodama}
Hideo Kodama. Quantum Gravity by the Complex Canonical Formulation.
gr-qc/9211022, Int.J.Mod.Phys.D1:439-524,1992.

\bibitem{Soo}
Chopin Soo. Wave function of the Universe and Chern-Simons Perturbation Theory. gr-qc/0109046, Class.Quant.Grav. 19 (2002) 1051-1064.

\bibitem{Pullin}
Jorge Pullin. Knot theory and quantum gravity in loop space: a primer. 
hep-th/9301028, AIP Conf.Proc.317:141-190,1994


\bibitem{WittKodama} E. Witten, A note on the Chern-Simons and Kodama wave functions,
arXiv:gr-qc/0306083.

\bibitem{Witten} E. Witten, Quantum field Theory and the Jones Polynomial, {\em
Commun. Math. Phys.},vol. 121, 1989, pp. 351-399.



\end{thebibliography}
\end{document}